\documentclass[12pt,a4paper]{amsart}
\usepackage{amsfonts,amssymb,amscd,amsmath,enumerate,verbatim,calc}
\usepackage{mathrsfs}
\usepackage[all]{xy}
\headsep=1cm \numberwithin{equation}{section}
\newtheorem{thm}{Theorem}[section]
\newtheorem{cor}[thm]{Corollary}
\newtheorem{lem}[thm]{Lemma}

\theoremstyle{definition}

\theoremstyle{remark}

\numberwithin{equation}{section}

\newcommand\Supp{\operatorname{Supp}}
\newcommand\Ass{\operatorname{Ass}}

\newcommand\mAss{\operatorname{mAss}}

\newcommand\ara{\operatorname{ara}}

\newcommand\cd{\operatorname{cd}}

\newcommand\Ext{\operatorname{Ext}}
\newcommand\Tor{\operatorname{Tor}}

\newcommand\height{\operatorname{height}}

\newcommand\Max{\operatorname{Max}}

\newcommand\m{\operatorname{\frak m}}

\newcommand\p{\operatorname{\frak p}}
\newcommand\q{\operatorname{\frak q}}

\begin{document}\title[Abelian categories of cofinite modules ]{A note on Abelian categories of cofinite modules }
\author[K. Bahmanpour ]{  Kamal Bahmanpour  }

\address{Department of Mathematics, Faculty of Sciences, University of Mohaghegh Ardabili,
56199-11367, Ardabil, Iran;
and School of Mathematics, Institute for Research in Fundamental Sciences (IPM), P.O. Box. 19395-5746, Tehran, Iran.} \email{\it bahmanpour.k@gmail.com}

\thanks{ 2010 {\it Mathematics Subject Classification}: Primary 13D45; Secondary 14B15, 13E05.\\This research of the author was supported by a grant from IPM (No. 96130018).}
\keywords{Abelian category, cofinite module, cohomological dimension, local
cohomology, Noetherian  ring.}

\begin{abstract}
   Let $R$ be a commutative Noetherian ring and $I$ be an ideal of $R$. In this article we answer affirmatively a question
raised by the present author in \cite{B2}. Also, as an immediate consequence of this result it is shown that the category of all $I$-cofinite $R$-modules $\mathscr{C}(R, I)_{cof}$  is an Abelian subcategory of the category of all $R$-modules, whenever $q(I,R)\leq 1$. These assertions answer affirmatively a question
raised by R. Hartshorne in [{\it Affine duality and cofiniteness},
Invent. Math. {\bf9}(1970), 145-164], in some special cases.
  \end{abstract}
\maketitle
\section{Introduction}

Throughout this article, let $R$ denote a commutative Noetherian ring
(with identity) and $I$ be an ideal of $R$. For an $R$-module $M$, the
$i$th local cohomology module of $M$ with support in $V(I)$
is defined as:
$$H^i_I(M) = \underset{n\geq1} {\varinjlim}\,\, \Ext^i_R(R/I^n,
M).$$   We refer the reader to \cite{BS} or \cite{Gr1} for more
details about local cohomology.\\

For an $R$-module $M$, the notion $\cd (I, M)$, {\it the
 cohomological dimension of $M$ with respect to $I$}, is defined as:
$$\cd(I,M)={\rm sup}\{i\in \Bbb{N}_0\,\,:\,\,H^i_I(M)\neq 0\}$$
and the notion $q(I,M)$, which for first time was introduced by Hartshorne, is defined as:
 $$q(I,M)={\rm sup}\{i\in \Bbb{N}_0\,\,:\,\,H^i_I(M)\,\,{\rm is }\,\,{\rm not }\,\,{\rm Artinian } \},$$  with the usual convention that the
supremum of the empty set of integers is interpreted as $-\infty$. These two notions have been studied by
several authors (see \cite{B, DY, DNT, Fa, GBA, Ha2, HL}). \\

Hartshorne in \cite{Ha}  defined an $R$-module $X$ to be
$I$-{\it cofinite}, if $\Supp X\subseteq
V(I)$ and ${\rm Ext}^{i}_{R}(R/I, X)$ is a finitely generated $R$-module
for each integer $i\geq0$. Then he posed the following question:\\

{\bf Question 1:} {\it Whether the category $\mathscr{C}(R, I)_{cof}$
of $I$-cofinite modules is an Abelian subcategory of the category of all $R$-modules?
That is, if $f: M\longrightarrow N$ is an $R$-homomorphism of
$I$-cofinite modules, are $\ker f$ and ${\rm coker} f$ $I$-cofinite}?\\

With respect to the question (1), Hartshorne gave a counterexample to show that this question has not an
affirmative answer in general, (see \cite[Section 3]{Ha}). On the positive side, Hartshorne proved that if $I$ is a prime ideal of dimension one in a complete regular local ring $R$, then the answer to his question
is yes. Delfino and Marley extended this result to arbitrary
complete local rings (see \cite{DM}). Kawasaki generalized the Delfino and Marley's result
for an arbitrary ideal $I$ of dimension one in a local ring $R$ (see \cite{Ka2}). Melkersson removed the local condition on the ring (see \cite{Me1}). Finally, in \cite{BNS} as a generalization of
 Melkersson's result it is shown that for any ideal $I$ in any Noetherian ring $R$, the category $\mathscr{C}^1(R, I)_{cof}$ of all $I$-cofinite $R$-modules $M$ with $\dim M\leq 1$ is Abelian. For some other similar results,  see also \cite{B3}.\\

Recall that, for any proper ideal $I$ of $R$, the
\emph{arithmetic rank} of $I$, denoted by $\ara(I)$, is the least number of
elements of $I$ required to generate an ideal which has the same radical as $I$.\\

 Kawasaki
 proved that if $\ara(I)=1$ then the category $\mathscr{C}(R, I)_{cof}$
 is Abelian (see \cite{Ka1}). Pirmohammadi et al. in \cite{PAB} as a generalization of Kawasaki's result proved that if $I$ is an ideal of a Noetherian local ring with $\cd(I,R)\leq1$, then $\mathscr{C}(R, I)_{cof}$ is Abelian.  Recently, Divaani-Aazar et al. in \cite{DFT} have removed
the local condition on the ring. Finally, the present author in \cite{B} proved that if $I$ is
an ideal of a  Noetherian complete local ring $R$ with $q(I,R)\leq1$, then $\mathscr{C}(R, I)_{cof}$ is Abelian.\\

We recall that the present author in \cite{B2}, for any ideal $I$ of $R$ and any finitely generated $R$-module $M$, defined: $$\mathfrak{A}(I,M):=\{\p \in \mAss_R M\,\,:\,\,I+\p=R\,\,\,\,{\rm or}\,\,\,\,\p\supseteq I\},$$$$\mathfrak{B}(I,M):=\{\p \in \mAss_R M\,\,:\,\,\cd(I,R/\p)=1\},$$$$\mathfrak{C}(I,M):=\{\p \in \mAss_R M\,\,:\,\,q(I,R/\p)=1\}\,\,{\rm and}$$$$\mathfrak{D}(I,M):=\{\p \in \mAss_R M\,\,:\,\,0\leq \dim R/(I+\p)\leq 1\}.$$

He proved that if $\mAss_R R=\mathfrak{A}(I,R)\cup\mathfrak{B}(I,R)\cup\mathfrak{D}(I,R)$ then  $ \mathscr{C}(R,I)_{cof} $ is Abelian.  Also, he asked the following question (see \cite[Question D]{B2}):\\

{\bf Question 2:} {\it Let $R$ be a Noetherian ring and $I$ be an ideal of $R$ such that $$\mAss_R R= \mathfrak{A}(I,R) \cup\mathfrak{B}(I,R)\cup\mathfrak{C}(I,R)  \cup\mathfrak{D}(I,R).$$ Whether $\mathscr{C}(R, I)_{cof}$ is Abelian}? \\

In this article, we present an affirmative answer to Question {\rm 2}.  More precisely, we prove the following theorem:\\

{\bf Theorem 1.} {\it Let $ R $ be a Noetherian ring and $I$ be an ideal of $R$ such that $$\mAss_R R=\mathfrak{A}(I,R)\cup\mathfrak{B}(I,R)\cup\mathfrak{C}(I,R)\cup\mathfrak{D}(I,R).$$
Then $ \mathscr{C}(R,I)_{cof} $ is Abelian.}\\

Also, as an immediate consequence of Theorem 1, we deduce the following generalization of \cite[Theorem 5.3]{B} and \cite[Theorem 2.2]{DFT}.\\

{\bf  Corollary. } {\it Let $ R $ be a Noetherian ring and $I$ be an ideal of $R$ with $q(I,R)\leq 1$.
Then $ \mathscr{C}(R,I)_{cof} $ is Abelian.}\\

Throughout this paper, for any ideal $I$ of a Noetherian ring $R$, we denote the category of all $I$-cofinite $R$-modules by $\mathscr{C}(R, I)_{cof}$. Also, for each $R$-module $L$, we denote by
  $\mAss_R L$, the set of
 minimal elements of $\Ass_R L$ with respect to inclusion. For any ideal $\frak{a}$ of $R$, we denote
$\{\frak p \in {\rm Spec}\,R:\, \frak p\supseteq \frak{a}\}$ by
$V(\frak{a})$. For any Noetherian local ring $(R,\m)$, we denote by $\widehat{R}$ the $\m$-adic completion of $R$. Finally, we denote by $\Max(R)$ the set of all maximal ideals of $R$. For any unexplained notation and terminology we refer the reader to \cite{BS} and \cite{Mat}.

\section{Results}

The main purpose of this section is to prove Theorem 2.9, which presents an affirmative answer to a question
raised by the present author in \cite{B2}. The following auxiliary lemmas are quite useful in the proof of  Theorem 3.9.

\begin{lem} $(${\rm See} \cite[Corollary 3.6]{B2}$)$
\label{2.1}
Let $(R,\m)$ be a Noetherian complete local ring and $I$ be an ideal of $R$ such that $$\mAss_R R=\mathfrak{A}(I,R)\cup\mathfrak{B}(I,R)\cup\mathfrak{C}(I,R)\cup\mathfrak{D}(I,R).$$  Then
$\mathscr{C}(R, I)_{cof}$ is Abelian.
\end{lem}
\qed\\

\begin{lem}
\label{2.2}
Let $ (R,\frak m) $ be a Noetherian local ring and $T$ be an $R$-module such that the $\widehat{R}$-module $T\otimes_R \widehat{R}$ is finitely generated. Then the $R$-module $T$ is finitely generated.
\end{lem}
 \proof The assertion easily follows from the fact that $\widehat{R}$ is a faithfully flat $R$-algebra.\qed\\

\begin{lem} $(${\rm See} \cite[Lemma 3.1]{B2}$)$
\label{2.3}
Let $R$ be a Noetherian ring, $I$ be an ideal of $R$ and $M$ be a finitely generated $R$-module. Then,
$\mathfrak{B}(I,M)\subseteq \big(\mathfrak{C}(I,M)\cup\mathfrak{D}(I,M)\big)$. In particular,
$$\mAss_R M=\mathfrak{A}(I,M) \cup \mathfrak{B}(I,M)\cup \mathfrak{C}(I,M)\cup \mathfrak{D}(I,M)$$ if and only if $\mAss_R M=\mathfrak{A}(I,M) \cup \mathfrak{C}(I,M) \cup \mathfrak{D}(I,M).$
\end{lem}\qed\\

\begin{lem}
\label{2.4}
Let $(R,\m)$ be a Noetherian local ring and $I$ be an ideal of $R$ such that $$\mAss_R R=\mathfrak{A}(I,R)\cup\mathfrak{B}(I,R)\cup\mathfrak{C}(I,R)\cup\mathfrak{D}(I,R).$$  Then,
$$\mAss_{\widehat{R}} \widehat{R}=\mathfrak{A}(I\widehat{R},\widehat{R})\cup\mathfrak{C}(I\widehat{R},\widehat{R})\cup\mathfrak{D}(I\widehat{R},\widehat{R}).$$
\end{lem}
\proof
Since, by the hypothesis we have $$\mAss_R R=\mathfrak{A}(I,R)\cup\mathfrak{B}(I,R)\cup\mathfrak{C}(I,R)\cup\mathfrak{D}(I,R),$$ it follows from
Lemma 2.3 that $$\mAss_R R=\mathfrak{A}(I,R)\cup\mathfrak{C}(I,R)\cup\mathfrak{D}(I,R).$$
Now, in order to prove the assertion, it is enough to prove that
$$\mAss_{\widehat{R}} \widehat{R}\subseteq\mathfrak{A}(I\widehat{R},\widehat{R})\cup\mathfrak{C}(I\widehat{R},\widehat{R})\cup\mathfrak{D}(I\widehat{R},\widehat{R}).$$
Assume that $\mathfrak{P}\in \mAss_{\widehat{R}} \widehat{R}$ and set $\q:=\mathfrak{P}\cap R$. Then, there exists $\p\in \mAss_R R $ such that $\p\subseteq \q$. From the hypothesis it follows that $$\p \in \big(\mathfrak{A}(I,R)\cup\mathfrak{C}(I,R)\cup\mathfrak{D}(I,R)\big).$$  We consider the following three cases:\\

Case 1. Assume that $\p \in \mathfrak{A}(I,R)$. If $I+\p=R$ then, $$\widehat{R}=I\widehat{R}+\p \widehat{R}\subseteq I\widehat{R}+\mathfrak{P} \subseteq \widehat{R}$$ and hence $I\widehat{R}+\mathfrak{P}=\widehat{R}$. Therefore, $\mathfrak{P}\in \mathfrak{A}(I\widehat{R},\widehat{R})$. Also, if $\p\supseteq I$ then $I\widehat{R}\subseteq \p\widehat{R}\subseteq \mathfrak{P}$ and hence $\mathfrak{P}\in \mathfrak{A}(I\widehat{R},\widehat{R})$.\\

Case 2. Assume that $\p \in \mathfrak{C}(I,R)$. Then, as $q(I,R/\p)=1$ it follows that $$q(I\widehat{R}, \widehat{R}/\p \widehat{R})=1.$$ So, as $\Supp \widehat{R}/\mathfrak{P}  \subseteq \Supp \widehat{R}/\p \widehat{R},$ it follows from \cite[Theorem 3.2]{DY} that   $$q(I\widehat{R},\widehat{R}/\mathfrak{P} )\leq q(I \widehat{R},\widehat{R}/\p \widehat{R})= 1.$$ If $q(I\widehat{R},\widehat{R}/\mathfrak{P} )=1$, then  $\mathfrak{P}\in \mathfrak{C}(I \widehat{R},\widehat{R})$. Also, if $q(I\widehat{R},\widehat{R}/\mathfrak{P})=0$ then it is clear that $\mathfrak{P}\in V(I \widehat{R})$ and hence $\mathfrak{P}\in \mathfrak{A}(I \widehat{R},\widehat{R})$. Finally, if $q(I \widehat{R},\widehat{R}/\mathfrak{P})=-\infty$ then in view of \cite[Lemma 4.1]{B} we have $\mathfrak{P}\in \mathfrak{D}(I \widehat{R},\widehat{R})\cup\mathfrak{A}(I \widehat{R},\widehat{R})$.\\

Case 3. Assume that $\p \in \mathfrak{D}(I,R)$. Then, by the definition we have $$\dim \widehat{R}/(I\widehat{R}+\mathfrak{P})\leq \dim \widehat{R}/(I \widehat{R}+\p \widehat{R})= \dim R/(I+\p)\leq 1,$$which implies that $\mathfrak{P}\in \mathfrak{D}(I \widehat{R},\widehat{R})\cup\mathfrak{A}(I\widehat{R},\widehat{R})$.
\qed\\

\begin{lem}
\label{2.5}
Let $R$ be a Noetherian ring and $I$ be an ideal of $R$ such that $$\mAss_R R=\mathfrak{A}(I,R)\cup\mathfrak{B}(I,R)\cup\mathfrak{C}(I,R)\cup\mathfrak{D}(I,R).$$  Then for each multiplicative subset $S$ of $R$ we have
$$\mAss_{S^{-1}R} S^{-1}R=\mathfrak{A}(S^{-1}I,S^{-1}R)\cup\mathfrak{B}(S^{-1}I,S^{-1}R)\cup\mathfrak{C}(S^{-1}I,S^{-1}R)\cup\mathfrak{D}
(S^{-1}I,S^{-1}R).$$
\end{lem}
\proof  The proof is straightforward and is left to the reader.\qed\\

\begin{lem}
\label{2.6}
Let $ (R,\frak m) $ be a Noetherian local ring and $ I $ be an ideal of $ R $ such that
$$\mAss_R R=\mathfrak{A}(I,R)\cup\mathfrak{B}(I,R)\cup\mathfrak{C}(I,R)\cup\mathfrak{D}(I,R).$$
Then $\mathscr{C}(R,I)_{cof}$ is Abelian.
\end{lem}
\proof Let $M,\,N\in \mathscr{C}(R, I)_{cof}$ and let $f:M\longrightarrow N$ be
an $R$-homomorphism. It is enough to prove that the $R$-modules
${\rm ker} f$ and ${\rm coker} f$ are $I$-cofinite. Set $K:={\rm ker} f$. Since, $$\mAss_R R=\mathfrak{A}(I,R)\cup\mathfrak{B}(I,R)\cup\mathfrak{C}(I,R)\cup\mathfrak{D}(I,R),$$ from Lemma 2.4 it follows that
$$\mAss_{\widehat{R}} \widehat{R}=\mathfrak{A}(I\widehat{R},\widehat{R})\cup\mathfrak{C}(I\widehat{R},\widehat{R})\cup\mathfrak{D}(I\widehat{R},\widehat{R}).$$ Moreover, from the hypothesis $M,\,N\in \mathscr{C}(R, I)_{cof}$ it follows that $$M\otimes_R \widehat{R},\,N\otimes_R \widehat{R}\in \mathscr{C}(\widehat{R}, I\widehat{R})_{cof}.$$ Since, $\widehat{R}$ is a flat $R$-algebra it follows that $K\otimes_R \widehat{R}={\rm ker} (f\otimes_R \widehat{R})$. Therefore, Lemma 2.1 yields that the $\widehat{R}$-module $K\otimes_R \widehat{R}$ is $I\widehat{R}$-cofinite. Hence, for each integer $j\geq0$ the $\widehat{R}$-module $$\Ext^j_{\widehat{R}}(\widehat{R}/I\widehat{R},K\otimes_R \widehat{R})\simeq \Ext^j_R(R/I,K)\otimes_R \widehat{R}$$ is finitely generated. Now, Lemma 2.2 implies that for each integer $j\geq0$ the $R$-module $\Ext^j_R(R/I,K)$ is finitely generated, which means that $K$ is $I$-cofinit. Now, the assertion follows from the exact sequences
$$0\longrightarrow {\rm ker} f \longrightarrow M \longrightarrow
{\rm im} f \longrightarrow 0,$$ and $$0\longrightarrow {\rm im}f \longrightarrow N \longrightarrow
{\rm coker}f \longrightarrow 0.$$
\qed\\

\begin{lem} $(${\rm See} \cite[Theorem 2.5]{B2}$)$
\label{2.7}
Let $ R $ be a Noetherian ring and $I$ be an ideal of $R$ such that $$\mAss_R R=\mathfrak{A}(I,R)\cup\mathfrak{B}(I,R)\cup\mathfrak{D}(I,R).$$  Then
$\mathscr{C}(R, I)_{cof}$ is Abelian.
\end{lem}
\qed\\

\begin{lem}
\label{2.8}
Let $ R $ be a Noetherian ring and $J=(x_1,\dots,x_k)$ be an ideal of $R$. Let $\p_1, \p_2,..., \p_n$ be prime ideals of $R$ such that $J\not\subseteq \bigcup_{i=1}^n \p_i$. Then, there are elements $y_1,\dots,y_k\in J$ such that $J=(y_1,\dots,y_k)$ and $\left(\bigcup_{i=1}^n \p_i\right)\bigcap\{y_1,\dots,y_k\}=\emptyset$.
\end{lem}
\proof  Since, $J=(x_1,\dots,x_k)\nsubseteq \bigcup_{i=1}^n \p_i$ it follows that there is $a_1\in (x_2,\dots,x_k)$ such that
$$x_1+a_1\not\in \bigcup_{i=1}^n \p_i.$$
Set $y_1:=x_1+a_1$. Then, $J=
(y_1,x_2,\dots,x_k)$. We shall construct the sequence
$y_1,\dots,y_k$  which are not belong to $\bigcup_{i=1}^n \p_i$ and $J=(y_1,\dots,y_k)$ by an inductive process. To do this
end, assume that $1\leq i<k$, and that we have already constructed
elements $y_1,\dots,y_i$ such that $J= (y_1,\dots, y_i,
x_{i+1},\dots,x_k).$ We show how to construct $y_{i+1}$.\\

 To do
this, as $$J= (y_1,\dots, y_i,
x_{i+1},\dots,x_k)\nsubseteq
\bigcup_{i=1}^n \p_i$$
it follows that there is $a_{i+1}\in (y_1,\dots, y_i,
x_{i+2},\dots,x_k)$ such that $$x_{i+1}+a_{i+1}\not\in
\bigcup_{i=1}^n \p_i.$$

Set $y_{i+1}:=x_{i+1}+a_{i+1}$. Then, $J=
(y_1,\dots, y_i, y_{i+1},x_{i+2} \dots,x_k).$ This completes the
inductive step in the construction.\qed\\

Now, we are ready to state and prove the main result of this article, which answers affirmatively \cite[Question D]{B2}.

\begin{thm}
\label{2.9}
 Let $ R $ be a Noetherian ring and $I$ be an ideal of $R$ such that $$\mAss_R R=\mathfrak{A}(I,R)\cup\mathfrak{B}(I,R)\cup\mathfrak{C}(I,R)\cup\mathfrak{D}(I,R).$$
Then $\mathscr{C}(R,I)_{cof}$ is Abelian.
\end{thm}
\proof
Let $M,\,N\in \mathscr{C}(R, I)_{cof}$ and let $f:M\longrightarrow N$ be
an $R$-homomorphism. By the proof of Lemma 2.6 it is enough to prove that the $R$-module
${\rm ker} f$ is $I$-cofinite. Set $K:={\rm ker} f$. \\

If $\mAss_R R=\mathfrak{A}(I,R)\cup\mathfrak{B}(I,R)\cup\mathfrak{D}(I,R),$ then the assertion holds by Lemma 2.7. So, we may assume that $\mAss_R R\neq \mathfrak{A}(I,R)\cup\mathfrak{B}(I,R)\cup\mathfrak{D}(I,R).$\\

  Set $\Phi:=\mAss_R R \backslash\big( \mathfrak{A}(I,R)\cup\mathfrak{B}(I,R)\cup\mathfrak{D}(I,R)\big).$ Then, it is clear that $\Phi\subseteq \mathfrak{C}(I,R)$. Now, set $B:=\bigoplus_{\p\in \Phi}R/\p$  and assume that the ideal $I$ is generated by $n$ elements. Then, in view of \cite[Theorem 3.3.1]{BS} we have ${\rm cd}(I,B)\leq n$. Since, the $R$-module $H^i_I(B)$ is Artinian for each integer $i\geq 2$ it follows that $$W:=\bigcup_{i=2}^{{\rm cd}(I,B)}\Supp H^i_I(B)$$is a finite subset of $\Max(R)$. Assume that $W=\{\m_1,\m_2,...,\m_t\}.$ Then, using the {\it Grothendieck's Vanishing Theorem} we can deduce that $\height \m_c\geq 2$, for each $1\leq c\leq t$. Therefore, we have $\height \big(\bigcap_{c=1}^t\m_i\big)\geq 2$ and hence $$\bigcap_{c=1}^t\m_i\not\subseteq \bigcup_{\p\in \mAss_R R}\p.$$ Assume that $$\bigcap_{c=1}^t\m_c=(x_1,\dots ,x_k).$$ Then, by Lemma 2.8 there are elements $y_1,\dots ,y_k\in \bigcap_{c=1}^t\m_c$ such that $$\bigcap_{c=1}^t\m_c=(y_1,\dots ,y_k)\,\,\,\,{\rm and}\,\,\,\,y_{\ell}\not\in \bigcup_{\p\in \mAss_R R}\p,\,\,{\rm for}\,\,{\rm each}\,\,1\leq \ell \leq k.$$

   In particular, $y_{\ell}$ is not a nilpotent element of $R$ and hence $S_{\ell}=\{1_{_{R}},y_{\ell},y_{\ell}^2,y_{\ell}^3,...\}$ is a multiplicative subset of $R$, for each $1\leq \ell \leq k$.\\

Now, for each integer $i\geq0$, set $U_i:=\Ext^i_R(R/I,K)$. Then, for each $1\leq c \leq t$, since by the hypothesis we have $$\mAss_R R=\mathfrak{A}(I,R)\cup\mathfrak{B}(I,R)\cup\mathfrak{C}(I,R)\cup\mathfrak{D}(I,R),$$
from Lemma 2.5 it follows  that
$$\mAss_{R_{\m_{c}}} R_{\m_{c}}=\mathfrak{A}(IR_{\m_{c}},R_{\m_{c}})\cup\mathfrak{B}(IR_{\m_{c}},R_{\m_{c}})
\cup\mathfrak{C}(IR_{\m_{c}},R_{\m_{c}})\cup\mathfrak{D}(IR_{\m_{c}},R_{\m_{c}}).$$

 Moreover, since $M,\,N\in \mathscr{C}(R, I)_{cof}$ and $f:M\longrightarrow N$ is an $R$-homomorphism it follows that, for each $1\leq c \leq t$, we have $M_{\m_{c}},\,N_{\m_{c}}\in \mathscr{C}(R_{\m_{c}}, IR_{\m_{c}})_{cof}$ and $f_{\m_{c}}:M_{\m_{c}}\longrightarrow N_{\m_{c}}$ is an $R_{\m_{c}}$-homomorphism with the kernel $K_{\m_{c}}$. So, in view of Lemma 2.6
  for each $1\leq c \leq t$ and each $i\geq0$, the $R_{\m_{c}}$-module, $(U_i)_{\m_{c}}\simeq \Ext^i_{R_{\m_{c}}}(R_{\m_{c}}/IR_{\m_{c}},K_{\m_{c}})$ is finitely generated.\\

On the other hand, for each $1\leq \ell \leq k$ and each $j\geq 2$, since $$\Supp H^j_I(B)\subseteq \{\m_1,\m_2,...,\m_t\}\subseteq V(Ry_{\ell})$$ it follows that $$H^j_{IR_{y_{\ell}}}(B_{y_{\ell}})\simeq\big(H^j_I(B)\big)_{y_{\ell}}=(S_{\ell})^{-1}H^j_I(B)=0,$$which means that
 $$\mAss_{R_{y_{\ell}}} R_{y_{\ell}}=\mathfrak{A}(IR_{y_{\ell}},R_{y_{\ell}})\cup\mathfrak{B}(IR_{y_{\ell}},
 R_{y_{\ell}})\cup\mathfrak{D}(IR_{y_{\ell}},R_{y_{\ell}}).$$

 Furthermore, since $M,\,N\in \mathscr{C}(R, I)_{cof}$ and $f:M\longrightarrow N$ is an $R$-homomorphism it follows that $M_{y_{\ell}},\,N_{y_{\ell}}\in \mathscr{C}(R_{y_{\ell}}, IR_{y_{\ell}})_{cof}$ and $(S_{\ell})^{-1}f:M_{y_{\ell}}\longrightarrow N_{y_{\ell}}$ is an $R_{y_{\ell}}$-homomorphism with the kernel $K_{y_{\ell}}$. Therefore, by Lemma 2.7,  the $R_{y_{\ell}}$-module
$$(U_i)_{y_{\ell}}\simeq Ext^i_{R_{y_{\ell}}}(R_{y_{\ell}}/IR_{y_{\ell}},K_{y_{\ell}})$$is finitely generated, for each $1\leq \ell \leq k$ and each $i\geq 0$. Thus, for each $1\leq \ell \leq k$ and each $i\geq 0$, there is a finitely generated submodule $G_{i,\ell}$ of the $R$-module $U_i$ such that $(U_i)_{y_{\ell}}=(G_{i,\ell})_{y_{\ell}}$ and so, $$\big(U_i/G_{i,\ell}\big)_{y_{\ell}}=0.$$

Now, for each $i\geq 0$, set $X_i:=G_{i,1}+G_{i,2}+\cdots+G_{i,k}$. Then, for each $i\geq 0$, the $R$-module $X_i$ is a finitely generated submodule of $U_i$ such that $\big(U_i/X_i\big)_{y_{\ell}}=0$, for each $1\leq \ell \leq k$. In particular, the $R$-module $U_i/X_i$ is $Ry_{\ell}$-torsion, for each $1\leq \ell \leq k$ and each $i\geq 0$. Consequently, the $R$-module $U_i/X_i$ is $(y_1,\dots ,y_k)$-torsion, for each $i\geq 0$. So, we have $$\Supp U_i/X_i \subseteq V(Ry_1+\cdots+Ry_k)=V(\cap_{c=1}^t\m_c)=\{\m_1,..,,\m_t\},\,\,\,{\rm for}\,\,{\rm each}\,\,i\geq0.$$
Since, for each $1 \leq c \leq t$ and each $i\geq 0$, the $R_{\m_{c}}$-module $\big( U_i/X_i\big)_{\m_{c}}$ is finitely generated, it follows that $U_i/X_i$ is a finitely generated $R$-module and hence $U_i$ is finitely generated too. This means that $K$ is an $I$-cofinite $R$-module, as required.
\qed\\

The following consequence of Theorem 2.9 is a generalization of \cite[Theorem 5.3]{B} and \cite[Theorem 2.2]{DFT}.

\begin{cor}
\label{2.10}
 Let $ R $ be a Noetherian ring and $I$ be an ideal of $R$ with $q(I,R)\leq 1$.
Then $\mathscr{C}(R,I)_{cof}$ is Abelian.
\end{cor}
\proof
Since, by the hypothesis we have $q(I,R)\leq 1$, from \cite[Theorem 3.2]{DY} it follows that $q(I,R/\p)\leq 1$, for each $\p\in \mAss_R R$. Now, let $\p\in \mAss_R R$. If $q(I,R/\p)=1$, then  $\p\in \mathfrak{C}(I,R)$. Also, if $q(I,R/\p)=0$ then it is clear that $\p\in V(I)$ and hence $\p\in \mathfrak{A}(I,R)$. Also, if $q(I,R/\p)=-\infty$ then in view of \cite[Lemma 4.1]{B} we have $\p\in \mathfrak{D}(I,R)\cup\mathfrak{A}(I,R)$. Therefore, we have$$\mAss_R R=\mathfrak{A}(I,R)\cup\mathfrak{C}(I,R)\cup\mathfrak{D}(I,R),$$and hence the assertion holds by Theorem 2.9.\qed\\

\begin{cor}
\label{2.11}
 Let $ R $ be a Noetherian ring of dimension at most $2$ and $I$ be an arbitrary ideal of $R$. Then $\mathscr{C}(R,I)_{cof}$ is Abelian.
\end{cor}
\proof
Let $d:=\dim R$. Then, by \cite[Proposition 5.1]{Me} the $R$-module $H^d_I(R)$ is Artinian and by {\it Grothendieck's Vanishing Theorem}, for each $i> d$ we have $H^i_I(R)=0$. Since, by the hypothesis we have $d\leq 2$ it is clear that $q(I,R)\leq 1$ and hence the assertion follows from Corollary 2.10.\qed\\

\begin{cor}
 \label{2.12}
Let $R$ be a Noetherian ring and $I$ be an ideal of $R$ such that
$$\mAss_R R=\mathfrak{A}(I,R)\cup\mathfrak{B}(I,R)\cup\mathfrak{C}(I,R)\cup\mathfrak{D}(I,R).$$
Let $$X^\bullet:
\cdots\longrightarrow X^i \stackrel{f^i} \longrightarrow X^{i+1}
\stackrel{f^{i+1}} \longrightarrow X^{i+2}\longrightarrow \cdots,$$ be a
complex such that  $X^i\in\mathscr{C}(R, I)_{cof}$ for all $i\in\Bbb{Z}$.
Then for each $i\in\Bbb{Z}$ the $i$-th cohomology module $H^i(X^\bullet)$ is in
$\mathscr{C}(R, I)_{cof}$.
\end{cor}
\proof The assertion follows from  Theorem 2.9.\qed\\

\begin{cor}
\label{2.13}
Let $R$ be a Noetherian ring and $I$ be an ideal of $R$ such that $q(I,R)\leq 1$.
Let $$X^\bullet:
\cdots\longrightarrow X^i \stackrel{f^i} \longrightarrow X^{i+1}
\stackrel{f^{i+1}} \longrightarrow X^{i+2}\longrightarrow \cdots,$$ be a
complex such that  $X^i\in\mathscr{C}(R, I)_{cof}$ for all $i\in\Bbb{Z}$.
Then for each $i\in\Bbb{Z}$ the $i$-th cohomology module $H^i(X^\bullet)$ is in
$\mathscr{C}(R, I)_{cof}$.
\end{cor}
\proof  Using Corollary 2.12, the assertion follows from the proof of Corollary 2.10.\qed\\

\begin{cor}
\label{2.14}
 Let $ R $ be a Noetherian ring of dimension at most $2$ and $I$ be an ideal of $R$. Let $$X^\bullet:
\cdots\longrightarrow X^i \stackrel{f^i} \longrightarrow X^{i+1}
\stackrel{f^{i+1}} \longrightarrow X^{i+2}\longrightarrow \cdots,$$ be a
complex such that  $X^i\in\mathscr{C}(R, I)_{cof}$ for all $i\in\Bbb{Z}$.
Then for each $i\in\Bbb{Z}$ the $i$-th cohomology module $H^i(X^\bullet)$ is in
$\mathscr{C}(R, I)_{cof}$.
\end{cor}
\proof Using Corollary 2.13, the assertion follows from the proof of Corollary 2.11.\qed\\

\begin{cor}
\label{2.15}
Let $R$ be a Noetherian ring and $I$ be an ideal of $R$ such that
$$\mAss_R R=\mathfrak{A}(I,R)\cup\mathfrak{B}(I,R)\cup\mathfrak{C}(I,R)\cup\mathfrak{D}(I,R).$$ Let
 $M$ be an $I$-cofinite $R$-module and $N$ be a finitely generated $R$-module. Then the $R$-modules
$\Tor_i^R(N,M)$ and $\Ext^i_R(N,M)$ are $I$-cofinite, for all integers $i\geq0$.
\end{cor}
\proof Since $N$ is finitely generated it follows that, $N$ has a
free resolution with finitely generated free $R$-modules. Now the
assertion follows using Corollary 2.12 and computing the $R$-modules $\Tor_i^R(N,M)$ and $\Ext^i_R(N,M)$, by this
free resolution. \qed\\

\begin{cor}
\label{2.16}
Let $R$ be a Noetherian ring and $I$ be an ideal of $R$ such that $q(I,R)\leq 1$. Let
 $M$ be an $I$-cofinite $R$-module and $N$ be a finitely generated $R$-module. Then the $R$-modules
$\Tor_i^R(N,M)$ and $\Ext^i_R(N,M)$ are $I$-cofinite, for all integers $i\geq0$.
\end{cor}
\proof Using Corollary 2.15, the assertion follows from the proof of Corollary 2.10.\qed\\

\begin{cor}
\label{2.17}
 Let $ R $ be a Noetherian ring of dimension at most $2$ and $I$ be an ideal of $R$. Let
 $M$ be an $I$-cofinite $R$-module and $N$ be a finitely generated $R$-module. Then the $R$-modules
$\Tor_i^R(N,M)$ and $\Ext^i_R(N,M)$ are $I$-cofinite, for all integers $i\geq0$.
\end{cor}
\proof Using Corollary 2.16, the assertion follows from the proof of Corollary 2.11.\qed\\

\begin{cor}
\label{2.18}
 Let $ R $ be a Noetherian ring and $I$ be an ideal of $R$ such that $$\mAss_R R=\mathfrak{A}(I,R)\cup\mathfrak{B}(I,R)\cup\mathfrak{C}(I,R)\cup\mathfrak{D}(I,R).$$
Let $M,\,N$ be two finitely generated $R$-modules. Then the $R$-modules
$\Tor_i^R(N,H^j_I(M))$ and $\Ext^i_R(N,H^j_I(M))$ are $I$-cofinite, for all integers $i\geq0$ and $j\geq 0$.
\end{cor}
\proof The assertion follows from \cite[Theorem 3.8]{B2} and Corollary 2.15.\qed\\

\begin{cor}
\label{2.19}
 Let $R$ be a Noetherian ring and $I$ be an ideal of $R$ such that $q(I,R)\leq 1$.
Let $M,\,N$ be two finitely generated $R$-modules. Then the $R$-modules
$\Tor_i^R(N,H^j_I(M))$ and $\Ext^i_R(N,H^j_I(M))$ are $I$-cofinite, for all integers $i\geq0$ and $j\geq 0$.
\end{cor}
\proof The assertion follows from \cite[Theorem 4.10]{B2} and Corollary 2.16.\qed\\

\begin{cor}
\label{2.20}
Let $ R $ be a Noetherian ring of dimension at most $2$ and $I$ be an ideal of $R$.
Let $M,\,N$ be two finitely generated $R$-modules. Then the $R$-modules
$\Tor_i^R(N,H^j_I(M))$ and $\Ext^i_R(N,H^j_I(M))$ are $I$-cofinite, for all integers $i\geq0$ and $j\geq 0$.
\end{cor}
\proof Using the proof of Corollary 2.11 we have $q(I,R)\leq 1$. So, the assertion follows from Corollary 2.19.\qed\\

\subsection*{Acknowledgements}
The author would like to thank to School of Mathematics, Institute for Research in Fundamental
Sciences (IPM) for its financial support.


\end{document}